# Solving the waste bin location problem with uncertain waste generation rate: a bi-objective robust optimization approach

Diego Rossit[1] and Jonathan Bard[2]

[1] Department of Engineering, INMABB, Universidad Nacional del Sur (UNS)-CONICET
1253 Alem Avenue, Bahía Blanca B8000CPB, Argentina
diego.rossit@uns.edu.ar

[2] Graduate Program in Operations Research & Industrial Engineering, University of Texas at Austin
204 E. Dean Keeton Street, Austin TX78712, United States
jbard@utexas.edu

Abstract

An efficient Municipal solid waste (MSW) system is critical to modern cities in order to enhance sustainability and livability of urban life. With this aim, the planning phase of the MSW system should be carefully addressed by decision makers. However, planning success is dependent on many sources of uncertainty that can affect key parameters of the system, e.g., the waste generation rate in an urban area. With this in mind, this paper contributes with a robust optimization model to design the network of collection points (i.e., location and storage capacity), which are the first points of contact with the MSW system. A central feature of the model is a bi-objective function that aims at simultaneously minimizing the network costs of collection points and the required collection frequency to gather the accumulated waste (as a proxy of the collection cost). The value of the model is demonstrated by comparing its solutions with those obtained from its deterministic counterpart over a set of realistic instances considering different scenarios defined by different waste generation rates. The results show that the robust model finds competitive solutions in almost all cases investigated. An additional benefit of the model is that it allows the user to explore trade-offs between the two objectives.

## 1. Introduction

Modern cities are increasingly dependent on an effective Municipal solid waste (MSW) system to maintain the quality of life and standard of living that most of us have come to expect. There are at least three reasons why efficient waste processing has become a top priority for many local governments. First, the correct handling of municipal waste can enhance the sustainability of a city by supporting its recycling, reuse and recovery programs. Second, households and business are highly sensitive to negative performance of the MSW system. When it is functioning smoothly, rarely does anyone notice, but when services begin to break down (e.g., when collection points do not have the required storage capacity or when collections do not stay on schedule) the social, political and environmental impacts can be severe (Rossit and Nesmachnow, 2022). Finally, the cost of the MSW system usually represents a sizable fraction of the local budget even if the system is well run (Das and Bhattacharyya, 2015). However, poor management practices can lead to non-negligible cost overruns. For example, when waste is dumped in unsuitable locations, it must eventually be removed often at a much greater expense than had it been disposed of properly to begin with.

The importance of efficient planning cannot be overstated so it is vital to point out that one aspect that can affect the quality of the planning process is the accuracy of the input data. In MSW management, much of the input data that is used for decision making can be affected by uncertainty (Singh, 2019). In the broader context of managerial problems, Rosenhead et al. (1972) identified three categories in which an optimization and/or decision-making procedure can take place according to the quality of the available information: certainty, risk, and uncertainty. In the context of certainty, it is assumed that all model parameter values are known and deterministic optimization models are applied. In the context of risk, the input of the optimization process is affected by randomness but there is enough information to estimate the distributions of the probabilistic parameters. Stochastic optimization models would be applied in this case. Finally, the decision-making process can take place in a context of uncertainty, in which there is not enough data available to estimate the probabilistic distribution of the random parameters that are input to the model. In these circumstances robust optimization models are generally applied. The objective of both the stochastic and robust optimization models is to find a solution that would perform relatively well under any possible realization of the random parameters (Snyder, 2006). Thus, the use of these methodologies can help to take into account uncertainty in the decision-making process.

In the case of the MSW system, one of the parameters that can be affected by random- ness is the waste generation rate, which is a key element to estimate the required processing capacity of the system, in general, and of the collection network, in particular. There are several alternatives for structuring the collection network such as door-to-door (or curbside) collection, collection points/community bins or self-delivered system for specific types of waste (Rossit and Nesmachnow, 2022). Community bins or collection points have certain advantages over the door-to-door approach. If collection points are well placed, they reduce the number of points to visit (i.e., from visiting every dwelling to visiting only the collection points). This reduces the traveled distance of the collection vehicle and, thus, the greenhouse emissions (Blazquez and Paredes-Belmar, 2020; Gilardino et al., 2017). In addition, considering that door-to-door collection in many cities is generally related to manually lifting bags or small bins (Battini et al., 2018), collection points can reduce the amount of effort of the personnel since they are generally associated with automatized collection. Finally, collection points can reduce the collection time spent loading the waste of bins since instead of collecting several small bins and bags, the truck has only to unload a larger bin. For example, Carlos et al. (2019) found in an empirical study that emptying two waste bins 1.1m3 required 85% more time than emptying a 2.4$m^3$ bin.

In light of these factors, the primary contribution of this paper is the presentation and evaluation of a robust bi-objective mixed-integer linear programming (MILP) model that addresses waste generation uncertainty when designing a collection network based on com- munity bins. The model involves both locating bin collection points and determining the required storage capacity to install in each (i.e., number and type of bins) in accordance with the uncertain amount of waste that is generated. As far as we are concerned, the optimization problem of designing a network collection point considering uncertain waste generation has not been previously addressed in the related literature. Furthermore, two optimization criteria are considered: (i) the minimization of the total investment cost, and (ii) the minimization of the number of required visits to empty the bins over the planning horizon (i.e., collection frequency) which is a proxy of posterior routing cost (Rossit et al., 2018). The model also includes a restriction for providing a minimum quality of service to all the generators in the area. Test results are provided using real scenarios associated with an Argentinean city that is trying to plan for different levels of waste generation. In summary, the main contributions of this work are twofold: it effectively addresses a real-world case study, providing insights that can be directly applied to practical waste management challenges. Additionally, it represents a pioneering effort in the related literature, being the first, to the best of our knowledge, to address the computationally complex waste bin location problem considering an uncertain waste generation rate by means of robust optimization. We consider that this dual contribution can help to the understanding and implementation of effective waste management strategies in the face of real-world complexities such as uncertain waste generation.

The paper is structured as follows. Section 1.1 reviews the most relevant literature followed by Section 2 that describes the robust mathematical formulation. Section 3 presents the computational experiments. Finally, conclusions are drawn in Section 4 and opportunities for future research are outlined.

*1.1. Related work*

A growing number of researchers have addressed the issue of locating collection points in urban areas (Rossit and Nesmachnow, 2022) but only a few have proposed methodologies that take into account some level of uncertainty. Kim and Lee (2013, 2015a, 2015b) considered the variation of waste generation of dwellings for different days within the planning horizon, allowing the model to assign users to different collection points on each day. In a second approach they begin with deterministic demand and perform sensitivity analysis to gain a better understanding of the modifications to the network that would be required if waste generation changes. Thus, they solve the waste bins location problem for common values of the waste generation rate (Letelier et al., 2022; Rathore et al., 2020; Rossit et al., 2020, 2018; Toutouh et al., 2018, 2020). Other articles have considered a simpler strategy where a safety factor is applied to expected bin usage. In this approach, bins are sized so that a specified proportion of their capacity will remain unused to absorb potential increases in the waste generation rate (Boskovic and Jovicic, 2015; Ferronato et al., 2020). Finally, Jammeli et al. (2021) although considered a fixed waste generation rate, they considered a stochastic normally distributed number of households within the surrounding area of each collection point and developed a transformed MILP formulation for network design. None of the previous articles have considered an uncertain waste generation in their models for designing the collection network. As far as we are concerned, this is the first article that deals with a probabilistic waste generation rate and proposes a robust model for this problem. Although there are not applications of robust models for the collection point location problem, there has been work that applied robust methods to other stages of the MSW logistic chain. Habibi et al. (2017) presented an optimization model to design a regional MSW network considering apart from the common economic and environmental aspects, the pollution affecting the area nearby the facilities (e.g., landfills, transfer station). In their model, uncertainty concerned the amount of recyclable and non-recyclable waste generated at each population center. Solutions were obtained with a combination of a robust and a two-stage stochastic approach proposed by Mulvey et al. (1995). In first stage of that approach, the authors minimized a weighted sum of the expected value of the second stage objective function plus the deviation of the second stage solution from the deterministic solution. Increasing the weight of the deviation component of the objective function, forces the optimization process to produce solutions that may present higher expected total costs but lower second-stage cost-deviations.

Tirkolaee et al. (2018) developed a robust periodic capacitated arc routing model for tackling the uncertainty surrounding waste generation. The aim of the model was to minimize the total variable (distance-related) plus fixed cost of the collection fleet; solutions were found with a commercial optimizer. In a second paper (Tirkolaee et al., 2019), the authors added an objective related to minimizing the longest traveled distance of any collection vehicle to reduce the makespan of the routing plan. This bi-objective problem was solved with the $\varepsilon$-constraint method coupled with an invasive weed-type heuristic. In a third paper Tirkolaee et al. (2020) addressed a location-allocation-inventory problem in which the aim was to design a strategic MSW network. The objective was to minimize the total cost of installing the collection bins and building the processing/disposal facilities, plus the operational, transportation, and penalty costs of uncollected waste (resulting from an excessive amount of waste due to the stochastic nature of its generation), and environmental costs related to emissions from the disposal process. A robust optimization approach was used to find solutions.

After an analysis of the related literature, we believe that improved solutions to the waste collection points location problem can be obtained by dealing with the underlying uncertainty in a more robust manner. In particular, there has been little work done on the reverse logistics chain associated with MSW in general that exploits robust methods. In particular, there is no article that deals with considering an uncertain waste generation rate when designing a collection network. The intent of this paper is to address this void.

## 2. The waste bin location problem with uncertain waste generation

In this section we discuss the problem formulation, the robust approach implemented to consider uncertainty and the resolution approach that is used for addressing the model.

### 1.2. Deterministic Bi-objective Problem Formulation

The waste bin location problem that is addressed in this article considers two different optimization criteria: to minimize the cost of the network of bins, which involves the cost of opening a collection point in an area and the cost of installing bins (Cavallin et al., 2020), and the minimization of the required collection frequency to empty the bins as a proxy of the posterior routing costs (Rossit et al., 2020). These are clearly conflicting objectives (Rossit et al., 2020). The larger is the storage capacity installed in the garbage accumulation points the larger is the investment cost. However, the collection frequency can be reduced since more waste can be accumulated. Given the distribution of waste generators, the model considers which collection points must be opened and with which capacity in order to optimize the aforementioned criteria. Waste generators are considered at an aggregated level (in this case per block) as it is usual in similar problems to make the problem computationally tractable (Cavallin et al., 2020; Rossit et al., 2018, 2020; Toutouh et al., 2020; Tralhão et al., 2010; Coutinho-Rodrigues et al., 2012) - which as an extension of the capacitated location problem is considered NP-hard (Cornuéjols et al., 1991)-. Additionally, a restriction to the quality of service (QoS) provided to the citizens is implemented. This requires that a waste generator is not assigned to a collection point that is beyond a certain threshold distance from its location.

Another feature of this model is that it considers source classified waste. This means that users are required to deposit different types of waste into separate bins at the collection points based on a specified classification. Then, the different waste fractions are collected separately by the collection vehicles. The purpose of source classification is to facilitate subsequent recycling and/or recovery activities in the following stages of the MSW system. The classification of MSW can be based on its source (domestic, commercial, institutional, and street sweeping), composition (wet, dry, organic, inert, glass, paper, metal, etc.), or level of hazard (toxic, infectious, radioactive, corrosive, flammable, etc.). In this study, we adopted a composition-based classification which involves categorizing it into two fractions: the mixed fraction and the recyclable fraction. This classification approach has been used in various related studies (Burnley et al., 2007; Cavallin et al., 2020; Rossit et al., 2020). In Bahía Blanca, the city of the case study, the current MSW system does not implement source classification of waste. Therefore, by considering only two fractions, we aim to facilitate a gradual transition of the population towards the new waste management system and foster the practice of waste separation at the source. The mixed fraction comprises organic waste (primarily food residues), diapers, textiles, wood, rubber, leather, cork, and similar types of waste. This fraction is typically directed to the landfill for final disposal. On the other hand, the recyclable fraction comprises metal, glass, paper, cardboard, and plastic. Proper sorting and treatment of this fraction allow for the recovery of various valuable resources. It is important to note that this study does not consider other waste types such as batteries, hazardous waste, construction materials, pathogenic residues, waste from street sweeping, and pruning remains. In Bahía Blanca, the public collection system does not handle these waste types due to city regulations. Instead, specialized private collection services are required for their proper management.

In the development of the model, we make use of the following notation.

*1.2.1. Sets*

$F$ = set of collection frequencies (frequencies profiles) of bins.

$H$ = set of types of bins.

$I$ = set of generation points that are also potential collection points in which (waste) bins can be installed.

$M$ = set of waste fractions.

### 1.2.2. Parameters

$Acc_f$ = maximum number of days between two consecutive visits with frequency $f$.

$Ar_h$ = required area of bin of type $h$.

$b_{im}$ = waste generation in volume at generation point $i$ of waste fraction $m$.

$Cap_h$ = capacity of bin of type $h$.

$co_h$ = cost of installing a bin of type $h$.

$Co_i$ = cost of conditioning a generation point $i$ to be a collection point in which bins are to be installed.

$di_{ij}$ = walking distance between generation points $i$ and $j$.

$dis^{max}$ = maximum allowable walking distance from a generation point (not used as a collection point) to the assigned collection point

$L$ = available area for installing bins at any generation point

### 1.2.3. Variables

$fr_{imf}$ = (binary) 1 if frequency $f$ was chosen for the collection point $i$ and type of waste $m$, 0 otherwise

$linV_{ijmf}$ = (continuous) representing the product of variables $x_{ij}$ and $fr_{imf}$

$q_i$ = (binary) 1 if the generation point $i$ is used as a collection point, 0 otherwise

$v_{hmi}$ = (integer) number of bins of type $h$ for waste fraction $m$ that are installed at collection point $i$

$x_{ij}$ = (binary) 1 if the user point $i$ is assigned to the collection point $j$, 0 otherwise

### 1.2.4. Waste collection model

$$Min \sum_{i \in I} \sum_{m \in M} \sum_{h \in H} Co_h v_{hmi} + \sum_{i \in I} Co_i q_i \tag{1a}$$

$$Max \sum_{i \in I} \sum_{m \in M} \sum_{f \in F} Acc_f fr_{imf} \tag{1b}$$

Subject to

$$\sum_{i \in I} x_{ij} = 1, \quad \forall j \in I \tag{1c}$$

$$\sum_{h \in H} Ar_h v_{hmi} \leq Lq_i, \quad \forall j \in I \tag{1d}$$

$$\sum_{j \in I} b_{jm} \left( \sum_{f \in F} Acc_f linV_{jimf} \right) \leq \sum_{h \in H} Cap_h v_{hmi}, \quad \forall i \in I, m \in M \tag{1e}$$

$$linV_{jimf} \leq x_{ij}, \quad \forall i, j \in I, m \in M, f \in F \tag{1f}$$

$$linV_{jimf} \leq fr_{jmf}, \quad \forall i, j \in I, m \in M, f \in F \tag{1g}$$

$$linV_{jimf} \geq x_{ij} + fr_{jmf} - 1, \quad \forall i, j \in I, m \in M, f \in F \tag{1h}$$

$$di_{ij} x_{ij} \leq dis^{max}, \quad \forall i, j \in I \tag{1i}$$

$$\sum_{j \in I} \sum_{f \in F} linV_{jimf} = 1, \qquad \forall\, i \in I, m \in M \tag{1j}$$

$$\sum_{f \in F} fr_{imf} = q_i, \qquad \forall\, i \in I, m \in M \tag{1k}$$

$$fr, q, x \in \mathcal{B}, v \in \mathcal{Z}_0^+, linV \geq 0 \tag{1l}$$

The first objective function (1a) is intended to minimize the cost of opening collection points and installing the different types of bins. The second objective function (1b) aims at maximizing the number of days that waste can remain uncollected in bins. Equation (1c) specifies that each generator point must be assigned to exactly one collection point. Constraints (1d) ensure that the area occupied by the bins at collection point $j$ is less than or equal to the maximum available space at the collection point. Constraints (1e) limit the total amount of waste $m$ processed by generators assigned to collection point $j$ to be no greater than the capacity of the installed bins for that type of waste at the collection point. Constraints (1f) - (1h) represent the linearization of the product of the variables $x_{ij}$ and $fr_{imf}$ using the auxiliary variable $linV_{jimf}$. Constraints (1i) guarantee that a generator is not assigned to a collection point that is farther than a threshold maximum distance, $dis^{max}$. This provides a minimum QoS to the users. Equation (1j) requires that every generator has to be assigned to a collection point that has a set frequency, while Equation (1k) only permits collection points to be assigned frequency if it is opened. Finally, variable definitions are provided in Constraints (1l).

*1.3. Consideration of uncertainty: robust transformation*

Robust optimization techniques in order to find a robust feasible solution that remains feasible within the realizations of the input data in the uncertainty set impose a related cost on the optimal solution based on its degree of feasibility (Tirkolaee et al., 2020). Thus, despite being relatively more expensive, the derived robust solution is more reliable and will allow the decision-maker to have a practical solution for several realizations of the uncertain parameters which are difficult to predict. The first robust approach was proposed by Soyster (1973) using a linear optimization model that provides the best feasible solution for all possible realizations of random input data. This approach tends to find "over-conservative" solutions, which means that in order to ensure the robustness of the solutions in most realizations of the random input data the results are often suboptimal. To consider less conservative solutions, Ben-Tal and Nemirovski (1999, 2000) and El Ghaoui et al. (1998) proposed other approaches which involves solving the robust counterparts of the nominal problem using a quadratic objective function. Although these models can better approximate some types of uncertainties without defaulting to over-conservative solutions, they have the disadvantage of requiring the solution of a non-linear optimization problem, which tends to be more difficult to solve than the linear model of Soyster (1973). As an intermediate strategy to tackle both over-conservatism and high computational effort, Bertsimas and Sim (2004) introduced a methodology to control the conservatism level of the solutions with a linear formulation. We have adopted their methodology for our work.

On this basis, to apply Bertsimas and Sim's robust methodology to Model (1) we consider that waste generation $b_{im}$ varies uniformly within a given range $\tilde{b}_{im} \epsilon [\bar{b}_{im} - \rho\bar{b}_{im}, \bar{b}_{im} + \rho\bar{b}_{im}]$. Although, uniform variation of the stochastic parameters is an assumption of Bertsimas and Sim's robust methodology, it is also a usual assumption in other works that addressed optimization problem in waste management (Tirkolaee et al., 2018, 2019). For simplicity, we perform the following substitution $\hat{b}_{im} = \rho\bar{b}_{im}$ hereafter. $\rho$ is referred to as the uncertainty level as it quantifies the maximum variation of $\tilde{b}_{im}$. With this consideration, we replace deterministic Constraint (1e) with the following equations:

$$\sum_{i \in I} \bar{b}_{im} \left( \sum_{f \in F} Acc_f linV_{jimf} \right) + z_{im}\Gamma_{im} + \sum_{j \in J_{im}} p_{jm} \leq \sum_{h \in H} cap_h v_{hmi}, \quad \forall\, i \in I, m \in M \tag{2a}$$

$$p_{jm} + z_{im} \geq \sum_{j \in J_{im}} \left( \hat{b}_{jm} \left( \sum_{f \in F} Acc_f y_{jimf} \right) \right), \quad \forall\, i,j \in I, m \in M \tag{2b}$$

$$y_{jimf} \leq linV_{jimf} \leq y_{jimf} \tag{2c}$$

$$z, p, y \geq 0 \tag{2d}$$

Where $z, p, y$ are continuous variables that appear in the robust transformation and $\Gamma_{im}$ is the conservatism level. The conservatism level is linked to the fact that it is unlikely that all the parameters that can potentially vary will actually do so simultaneously. Thus, the conservatism level fixes the number of parameters that are allowed to vary simultaneously. For a detailed step-by-step explanation of the transformation of the Model (1) to a robust optimization model the interested reader can refer to Appendix A.

Hereafter, we consider the *Deterministic* Model as the model composed by Equations (1a)-(1l) and the *Robust* Model as the model composed by Equations (1a)-(1d), (1f)-(1l), and (2a)-(2d).

### 1.4. Bi-objective approach

For our computations, we used the ε-constraint method, which has been successfully applied to other deterministic (Rossit et al., 2017, 2020; Tralhão et al., 2010; Coutinho-Rodrigues et al., 2012; Mavrotas et al., 2015) and robust problems in MSW management (Tirkolaee et al., 2019), to solve both Deterministic and Robust Models. Instead of the original version of the ε-constraint developed by Haimes (1971), we applied the augmented ε-constraint method developed by Mavrotas and Florios (2013) which improves the original version since it ensures that the solution found is Pareto optimal when it converges.

## 3 Computational experiments

This Section presents the implementation details, a description of the problem instances, the computations, and an analysis of the results.

*1.5. Implementation details*

The model was implemented in Python using Pyomo as the modelling language and Gurobi v10 as the integer programming solver. The instances were solved on a Dell Powerdge- T440 server with 28 cores and 128 GB of RAM memory. All runs were performed using a time limit of 1200 sec for the solver. For handling geographic information, we used QGIS (QGIS Development Team, 2023).

The set of instances was built using real data based on a neighborhood of Bahía Blanca. The complete instance considered in this study consists of 190 waste generation points, which are potential locations for collection points. Waste generation of each generation point was estimated through detailed field work, taking into account both the population distribution and the historical waste generated by small shops and institutions such as universities (Cavallin et al., 2020). Large shops, office buildings, and businesses are required to hire their own private waste collection firms according to the city regulations and, thus, were not included in the analysis. The collection frequencies considered for emptying the bins are: every day ($Acc_f = 1$), every two days ($Acc_f = 2$), and every three days ($Acc_f = 3$) during the week. The walking distances were estimated using Open Source Routing Machine2[1]. For this study, bins with back loading were chosen, as the current collection company in the city has compactor trucks that can be adapted for collecting this type of bins. The considered bins were of two types: CTB 1100 and "rear loading bin", which occupy 1.38×1.00 m2 and 1.00×1.00 m2, respectively. The volume and prices are 1.1 and 1 m3, and US$493 and US$528, respectively. These prices include the maintenance cost and the prorated installment cost based on the fact that the bins are expected to have a useful life of ten years. The cost of opening a collection point is set to the cost of one bin of type CTB 1100.

---
[1] http://project-osrm.org/

To test the scalability of the model, instances of varying sizes were generated by selecting contiguous areas that are subsets of the full area with 190 waste generation points. Table 1 provides an overview of the instances and includes the number of constraints and variables in the optimization model when applied to each instance for both the deterministic and robust models.

**Table 1**: Description of instances.

| Instance | Number of waste generation points | Waste generation per collection point | | | | Number of constraints Det/robus | Number of variables Det/robus |
|---|---|---|---|---|---|---|---|
| | | Recyclable waste [m$^3$] | | Mixed waste [m$^3$] | | | |
| | | Avg. | Std. Dev. | Avg. | Std. Dev. | | |
| i16 | 16 | 0.554 | 0.34 | 0.143 | 0.072 | 4993/8577 | 1968/3568 |
| i58 | 58 | 0.556 | 0.465 | 0.227 | 0.121 | 64381/111477 | 24186/44602 |
| i73 | 73 | 0.783 | 0.467 | 0.208 | 0.123 | 101836/176442 | 38106/70372 |
| i126 | 126 | 0.656 | 0.436 | 0.178 | 0.117 | 302653/524917 | 112518/208278 |
| i190 | 190 | 0.341 | 0.359 | 0.149 | 0.112 | 687421/1192821 | 254790/472150 |

For each instance, a set of scenarios with different conservatism and uncertainty levels (i.e., $\Gamma$ and $\rho$, respectively) are used. The considered values for $\Gamma$ are 0.05, 0.1, 0.2, 0.3, 0.4, and 0.1, 0.2, 0.3, 0.4 for $\rho$.

## 1.6. *Computational tests for obtaining multi-objective solutions*

In this section we present the results of the application of augmented ε-constraint method over the set of instances to obtain feasible solutions. For each instance ten runs are performed with the aim of obtaining ten multi-objective solutions that explore the trade-off among objectives. For brevity, detailed results of the computational experimentation for obtaining approximation of the worst values of objectives within the Pareto frontier are presented in Appendix B.

For analyzing the computational effort, we present Table 2 which presents for each in- stance and each scenario: the number of the multi-objective solutions found and the average optimality gap estimated by Gurobi. The average gap was calculated based on all the multi- objective runs performed by Gurobi for each particular instance (this excludes the runs for obtaining the best and worst values of each objective within the Pareto frontier).

From Table 2, it can be observed that the size of the instance has a significant impact on the number of multi-objective solutions found and the average gap. Generally, as the instance size increases, the number of solutions found decreases, and the average gap increases. This trend is expected as larger instances tend to have a higher complexity and a larger search space (as can be depicted in the number of variables and constraints in Table 1), making it more challenging to find a diverse set of efficient solutions.

Regarding the effect of $\Gamma$ and $\rho$, increasing these parameters can also impact the number of multi-objective solutions found and the average gap. The increment in $\Gamma$ and $\rho$ implies a more restricted problem, where the model needs to accommodate a larger expected amount of waste in the same number of potential garbage accumulation points. This increased restriction can make it more difficult to find a variety of feasible solutions, resulting in a reduced number of solutions and potentially a larger average gap.

It is worth noting that in scenarios where the number of obtained multi-objective solutions is substantially reduced, the average gap may be less representative due to the limited number of solutions available for calculation.

As mentioned, despite being relatively more expensive (in terms of values of the objective functions), the derived robust solution is more reliable and will allow the decision-maker to have a practical solution for several realizations of the uncertain parameters which are difficult to predict. To assess this deterioration in the objective functions we present a comparison between the solutions of the robust and deterministic models. In summarizing the computations, we compare three representative solutions in Table 3. This table presents average values of the five instances. First, we compare the solution with the minimum cost of each robust scenario (each combination of $\rho$ and $\Gamma$) with the solution with the minimum cost of the deterministic scenario. Second, we compare the best solution in terms of collection frequency among the robust and deterministic models. Finally, we compare the best compromise solution of the robust scenarios against the ideal vector (or solution) for the deterministic case. The ideal vector is constructed by combining the best values achieved in any run for each objective function, which means that it is inevitably an unattainable point (Rossit et al., 2022). The procedure of selecting this best compromise solution among the robust solutions is depicted in Fig. 1. Each of the three comparisons shown gives the percentage relative distance of each particular objective among the solutions, i.e., $\delta_c$ for the cost and $\delta_f$ for the collection frequency, and the overall distance $\Delta$ considering the Euclidean norm expressed by Equation (3):

$$\Delta = \sqrt{\sum_{o \in O} \left( \frac{RobSol_o - DetSol_o}{DetSol_o} \right)} \qquad (3)$$

where $O$ is the set of objectives and $DetSol_o$ and $RobSol_o$ are the values of the deterministic and robust solutions for objective $o \in O$.

In Table 3, the comparison of solutions with minimal cost shows that for a given $\Gamma$, the larger the uncertainty level $\rho$ the larger is the distance from the robust solution to the deterministic solution in terms of cost. Similarly, for a given $\rho$, the larger $\Gamma$ the larger the percentage distance of the cost of the robust solutions in comparison to the deterministic solutions. The scenarios with larger $\Gamma$ and $\rho$ are more demanding since the conservatism level $\Gamma$ sets the number of points that are expected to vary their waste production rate and the uncertainty level $\rho$ controls the maximum variation that is considered for each individual point. Thus, when more waste is expected to be generated more collection points have to be opened and more bins have to be installed and maintained increasing the network cost. In this minimum cost solution, usually $\Delta$ and $\delta_f$ follow the same rule for less demanding instances and they increase when $\Gamma$ and $\rho$ also increase. the exception are the instances with large $\Gamma$ and/or $\rho$. This might be related to the inability of Gurobi to find good solutions for the more challenging scenarios; in those cases, the solutions are likely to be suboptimal.

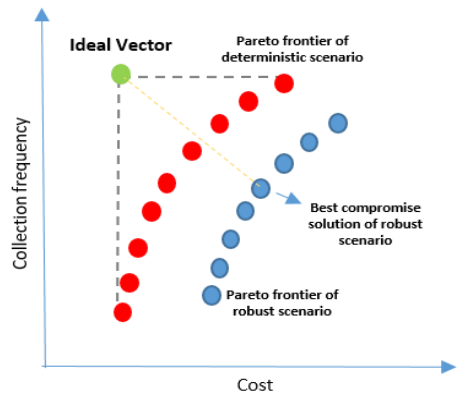

**Figure 1**: Best compromise solutions for the robust Pareto frontier.

For the comparisons among the solutions with maximum collection frequency, the behavior is similar. The larger $\Gamma$ and $\rho$, the larger is the distance between the robust and the deterministic solutions in terms of collection frequency. Similarly, the total expected generated waste of the scenario is larger with the increment of these parameters and this is handled by increasing the frequency with which collection points are emptied. However, the percentage distances of collection frequency are relatively smaller compared to the percentage distance of the cost in the previous minimum cost solutions. Excluding the particular case of the scenario with $\Gamma = 0.4$ and $\rho = 0.4$ (which has a $\delta_f$ of 32.63%), the largest $\delta_f$ is 1.02% while in the minimum cost solutions the largest $\delta_c$ is 41.31% (and the smallest $\delta_c$ is 3.69%). This may be connected to the fact that solutions with maximum frequency are relatively easier to obtain by the solver (as was evidenced in Appendix A) and the solver is able to find relatively good solutions for the robust scenarios.

When comparing the best compromise robust solution with the deterministic ideal vector (as is depicted in Figure 1) a similar pattern is observed with the overall distance Δ: the larger Γ and ρ, the greater the value of Δ. An important feature to highlight is that the ranges of variation in Δ within the same ρ or the same Γ are smaller in this compromise solution compared to the previous two extreme solutions. For example, in the compromise solution, the maximum range of variation for Δ occurs for Γ = 0.4, where Δ ranges from 92.26% (for ρ = 0.1) to 103.72% (for ρ = 0.4), resulting in a range of variation of 11.46%. In contrast, in the minimum cost solution, the largest range of variation for Δ is 56.34% (also when ρ = 0.4) , ranging from 27.29% (for Γ = 0.05) to 57.76% (for Γ = 0.4). Similarly, in the maximum collection frequency solution, the largest range of variation for Δ is also 56.34% (when ρ = 0.4), ranging from 7.85% (for Γ = 0.05) to 64.19% (for Γ = 0.4). This indicates that the compromise solution is relatively stable regarding this summary metric that considers both objectives simultaneously, which is a desirable characteristic for a compromise solution, when the scenario becomes more demanding due to a larger expected waste generation.

**Table 2**: Computational effort: Gurobi average gaps and number of multi-objective solutions found for each scenario.

| Instance | ρ = 0.1 | | | | | | | | | |
| --- | --- | --- | --- | --- | --- | --- | --- | --- | --- | --- |
| | Γ = 0.05 | | Γ = 0.1 | | Γ = 0.2 | | Γ = 0.3 | | Γ = 0.4 | |
| | Avg. gap | Sol | Avg. gap | Sol | Avg. gap | Sol | Avg. gap | Sol | Avg. gap | Sol |
| i16 | 13.59% | 10 | 14.93% | 10 | 16.82% | 10 | 17.94% | 10 | 18.69% | 10 |
| i58 | 18.44% | 10 | 19.37% | 10 | 21.94% | 10 | 23.70% | 10 | 24.60% | 10 |
| i73 | 22.33% | 10 | 22.37% | 10 | 24.81% | 10 | 26.65% | 10 | 28.30% | 10 |
| i126 | 25.03% | 10 | 25.58% | 10 | 27.21% | 8 | 29.75% | 9 | 30.19% | 8 |
| i190 | 31.28% | 9 | 33.25% | 8 | 35.20% | 8 | 30.29% | 5 | 33.19% | 5 |
| Average | 22.13% | 9.8 | 23.10% | 9.6 | 25.20% | 9.2 | 25.67% | 8.8 | 27.00% | 8.6 |
| | ρ = 0.2 | | | | | | | | | |
| Instance | Γ = 0.05 | | Γ = 0.1 | | Γ = 0.2 | | Γ = 0.3 | | Γ = 0.4 | |
| | Avg. gap | Sol | Avg. gap | Sol | Avg. gap | Sol | Avg. gap | Sol | Avg. gap | Sol |
| i16 | 13.94% | 10 | 16.10% | 10 | 18.17% | 10 | 18.72% | 10 | 19.47% | 10 |
| i58 | 18.87% | 10 | 21.02% | 10 | 26.28% | 10 | 28.65% | 10 | 29.80% | 8 |

| Instance | | | | | | | | | | |
|---|---|---|---|---|---|---|---|---|---|---|
| i73 | 22.77% | 10 | 24.54% | 10 | 28.40% | 10 | 30.74% | 10 | 28.75% | 5 |
| i126 | 26.18% | 10 | 27.75% | 10 | 30.34% | 8 | 32.50% | 7 | - | 0 |
| i190 | 23.13% | 3 | 32.96% | 6 | 23.64% | 2 | 23.81% | 1 | - | 0 |
| Average | 20.98% | 8.6 | 24.47% | 9.2 | 25.37% | 8 | 26.89% | 7.6 | 26.00% | 4.6 |

| | $\rho = 0.3$ | | | | | | | | | |
|---|---|---|---|---|---|---|---|---|---|---|
| Instance | $\Gamma = 0.05$ | | $\Gamma = 0.1$ | | $\Gamma = 0.2$ | | $\Gamma = 0.3$ | | $\Gamma = 0.4$ | |
| | Avg. gap | Sol | Avg. gap | Sol | Avg. gap | Sol | Avg. gap | Sol | Avg. gap | Sol |
| i16 | 15.40% | 10 | 16.55% | 10 | 18.28% | 10 | 20.21% | 10 | 20.71% | 10 |
| i58 | 20.73% | 10 | 24.01% | 10 | 28.81% | 10 | 29.96% | 7 | 30.11% | 7 |
| i73 | 23.64% | 10 | 26.09% | 10 | 30.63% | 9 | 28.44% | 4 | 27.29% | 3 |
| i126 | 28.10% | 10 | 29.22% | 9 | 28.38% | 4 | 30.31% | 3 | - | 0 |
| i190 | 31.54% | 6 | 31.52% | 5 | 26.38% | 1 | - | 0 | - | 0 |
| Average | 23.88% | 9.2 | 25.48% | 8.8 | 26.50% | 6.8 | 27.23% | 4.8 | 26.04% | 4.00 |

| | $\rho = 0.4$ | | | | | | | | | |
|---|---|---|---|---|---|---|---|---|---|---|
| Instance | $\Gamma = 0.05$ | | $\Gamma = 0.1$ | | $\Gamma = 0.2$ | | $\Gamma = 0.3$ | | $\Gamma = 0.4$ | |
| | Avg. gap | Sol | Avg. gap | Sol | Avg. gap | Sol | Avg. gap | Sol | Avg. gap | Sol |
| i16 | 15.50% | 10 | 17.30% | 10 | 18.82% | 10 | 20.54% | 10 | 21.94% | 10 |
| i58 | 21.47% | 10 | 26.14% | 10 | 32.13% | 10 | 26.44% | 4 | 33.75% | 5 |
| i73 | 24.66% | 10 | 28.43% | 10 | 32.49% | 8 | 33.25% | 6 | - | 0 |
| i126 | 27.87% | 10 | 30.30% | 8 | 23.41% | 1 | - | 0 | - | 0 |
| i190 | 31.65% | 6 | 49.92% | 2 | 29.20% | 2 | - | 0 | - | 0 |
| Average | 24.23% | 9.2 | 30.42% | 8 | 27.21% | 6.2 | 26.74% | 4 | 27.84% | 3 |

**Table 3**: Distance of representative solutions between the robust and the deterministic solutions.

| | Robust vs. deterministic solution with minimum cost | | | | | | | | | | | |
|---|---|---|---|---|---|---|---|---|---|---|---|---|
| $\Gamma$ | $\rho = 0.1$ | | | $\rho = 0.2$ | | | $\rho = 0.3$ | | | $\rho = 0.4$ | | |
| | $\Delta$ | $\delta_c$ | $\delta_f$ | $\Delta$ | $\delta_c$ | $\delta_f$ | $\Delta$ | $\delta_c$ | $\delta_f$ | $\Delta$ | $\delta_c$ | $\delta_f$ |
| 0.05 | 23.12% | 3.69% | 22.80% | 22.21% | 5.38% | 21.51% | 32.00% | 8.18% | 30.92% | 27.29% | 8.70% | 25.86% |
| 0.1 | 19.59% | 4.84% | 18.95% | 33.00% | 9.85% | 31.49% | 42.09% | 13.88% | 39.69% | 36.47% | 15.76% | 32.83% |
| 0.2 | 29.81% | 9.42% | 28.26% | 38.23% | 16.46% | 34.46% | 42.94% | 22.10% | 36.80% | 48.51% | 27.95% | 39.37% |
| 0.3 | 29.08% | 11.51% | 26.63% | 46.57% | 24.38% | 39.64% | 55.99% | 31.62% | 45.92% | 57.76% | 39.53% | 41.16% |
| 0.4 | 29.77% | 15.48% | 24.11% | 38.35% | 24.65% | 26.63% | 50.31% | 41.31% | 22.59% | 47.85% | 40.00% | 26.14% |
| | Robust vs. deterministic solution with maximum collection frequency | | | | | | | | | | | |
| $\Gamma$ | $\rho = 0.1$ | | | $\rho = 0.2$ | | | $\rho = 0.3$ | | | $\rho = 0.4$ | | |
| | $\Delta$ | $\delta_c$ | $\delta_f$ | $\Delta$ | $\delta_c$ | $\delta_f$ | $\Delta$ | $\delta_c$ | $\delta_f$ | $\Delta$ | $\delta_c$ | $\delta_f$ |
| 0.05 | 2.63% | 2.63% | 0.00% | 3.63% | 3.63% | 0.00% | 5.54% | 5.54% | 0.00% | 7.85% | 7.85% | 0.00% |
| 0.1 | 3.68% | 3.68% | 0.00% | 8.01% | 8.01% | 0.00% | 12.07% | 12.07% | 0.00% | 18.61% | 18.61% | 0.00% |
| 0.2 | 8.59% | 8.59% | 0.00% | 18.59% | 18.59% | 0.00% | 28.65% | 28.65% | 0.14% | 45.51% | 45.51% | 0.80% |
| 0.3 | 12.84% | 12.84% | 0.00% | 31.71% | 31.71% | 0.23% | 45.38% | 45.37% | 1.02% | 45.38% | 45.37% | 1.02% |
| 0.4 | 18.96% | 18.96% | 0.00% | 39.02% | 39.01% | 0.57% | 54.73% | 54.69% | 1.76% | 64.19% | 44.72% | 32.63% |
| | Robust compromising solution vs. deterministic ideal vector | | | | | | | | | | | |
| $\Gamma$ | $\rho = 0.1$ | | | $\rho = 0.2$ | | | $\rho = 0.3$ | | | $\rho = 0.4$ | | |
| | $\Delta$ | $\delta_c$ | $\delta_f$ | $\Delta$ | $\delta_c$ | $\delta_f$ | $\Delta$ | $\delta_c$ | $\delta_f$ | $\Delta$ | $\delta_c$ | $\delta_f$ |
| 0.05 | 87.63% | 36.62% | 79.01% | 85.98% | 36.86% | 76.74% | 86.58% | 37.16% | 77.69% | 87.33% | 37.89% | 78.17% |

| | | | | | | | | | | | |
|---|---|---|---|---|---|---|---|---|---|---|---|
| 0.1 | 88.34% | 37.82% | 79.20% | 89.64% | 39.59% | 78.04% | 90.60% | 48.61% | 75.46% | 98.56% | 61.80% | 74.05% |
| 0.2 | 91.87% | 50.00% | 76.72% | 93.27% | 52.81% | 75.92% | 97.37% | 59.69% | 75.87% | 98.87% | 56.52% | 79.49% |
| 0.3 | 94.03% | 53.56% | 76.80% | 98.05% | 62.65% | 73.36% | 100.85% | 67.21% | 73.49% | 93.99% | 49.64% | 79.55% |
| 0.4 | 96.18% | 56.70% | 76.98% | 92.04% | 52.81% | 74.99% | 97.90% | 62.68% | 74.72% | 107.14% | 69.81% | 81.06% |

### 1.7. Pareto frontiers of representative instances

To perform a visual study of Pareto frontiers of some illustrative instances we present in first place Figures 2-3 which are examples of the Pareto frontiers computed for a set of scenarios in instance i16 and i.58 respectively. It can be depicted that given an uncertainty level $\rho$ the successive frontiers with incremental $\Gamma$ tend to achieve worse values in terms of both objectives. In Figure 3d with $\rho = 0.4$ it becomes visible the reduction in the number of multi-objective solutions for large values of $\Gamma$ in demanding scenarios.

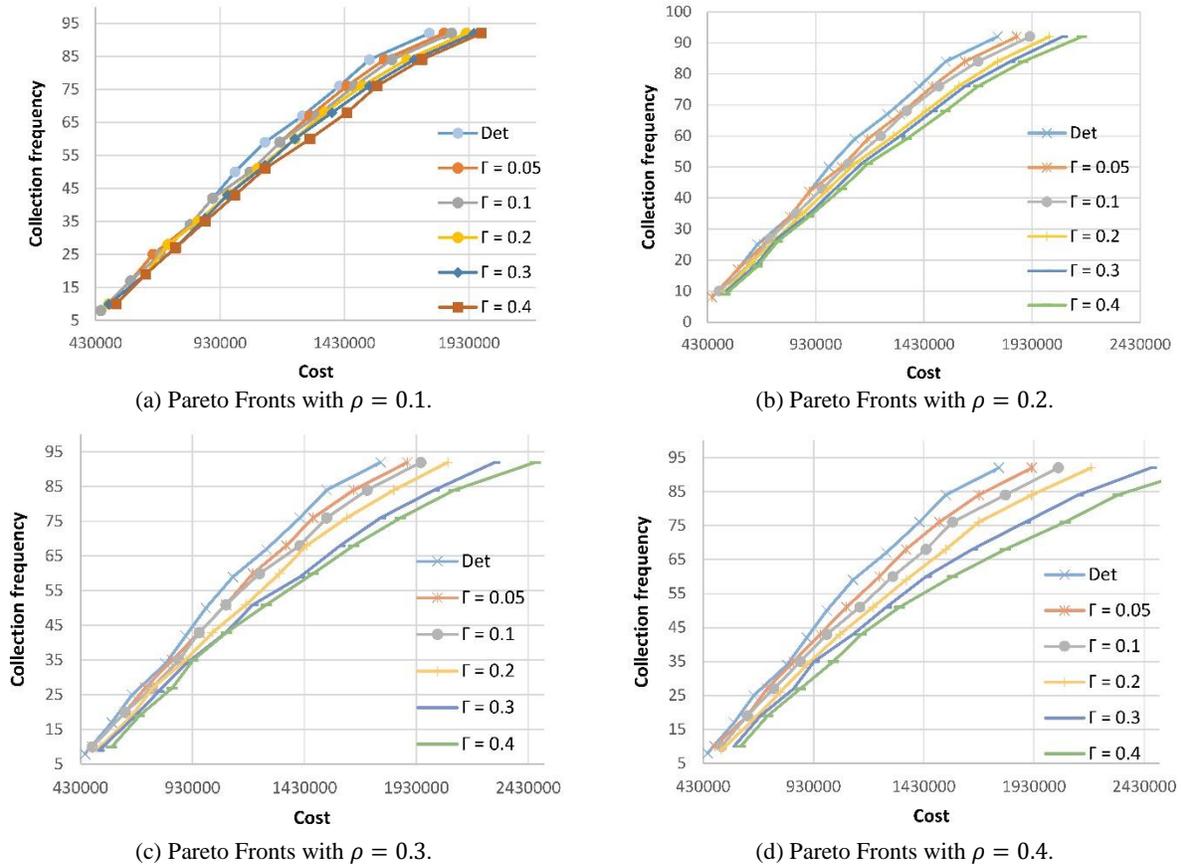

(a) Pareto Fronts with $\rho = 0.1$.
(b) Pareto Fronts with $\rho = 0.2$.
(c) Pareto Fronts with $\rho = 0.3$.
(d) Pareto Fronts with $\rho = 0.4$.

**Figure 2**: Solutions of instance i16.

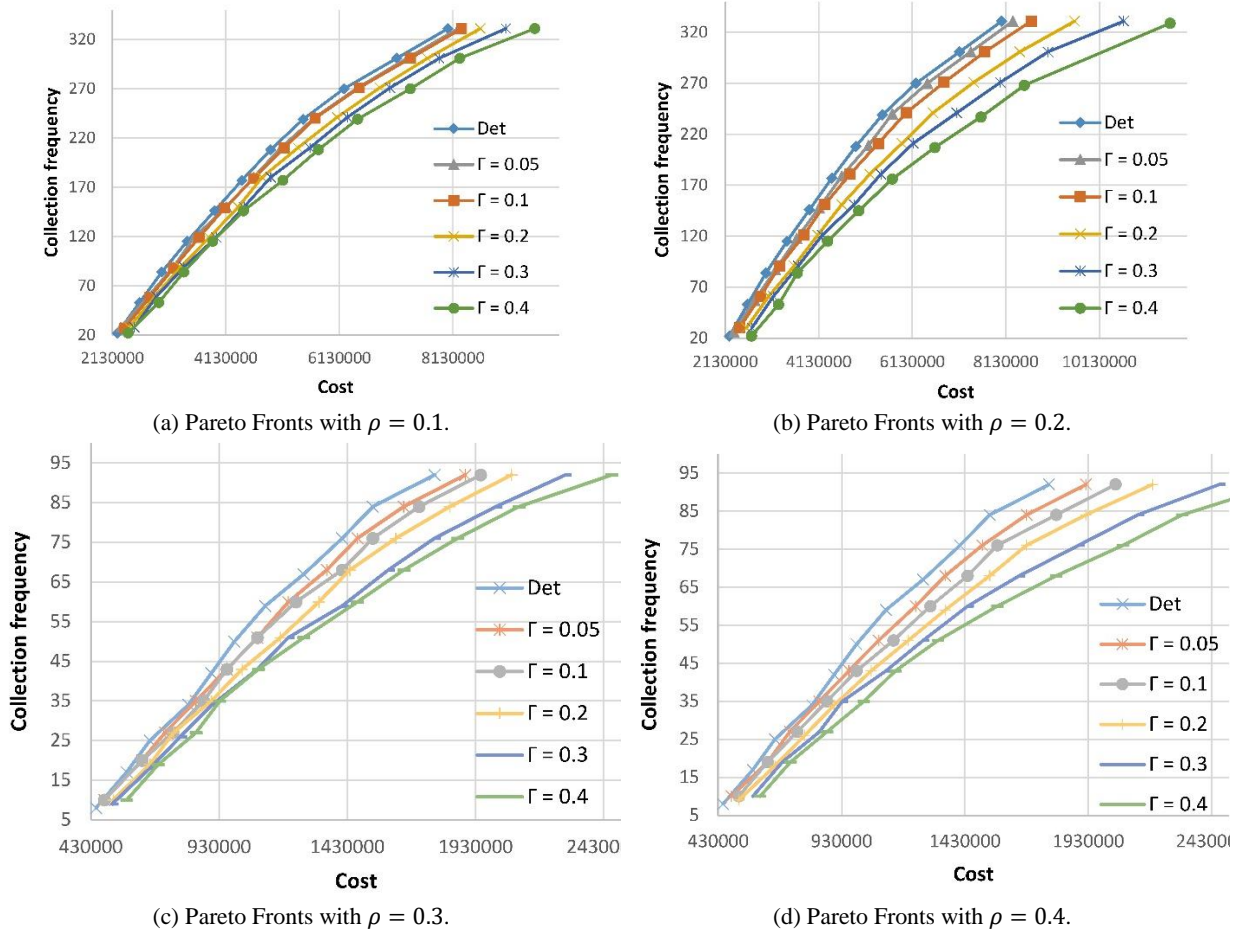

(a) Pareto Fronts with $\rho = 0.1$.
(b) Pareto Fronts with $\rho = 0.2$.
(c) Pareto Fronts with $\rho = 0.3$.
(d) Pareto Fronts with $\rho = 0.4$.

**Figure 3**: Solutions of instance i58.

*1.8. Geographic visualization of representative solutions*

Figure 4 presents examples of geographic visualizations of the solutions for instance i73 on the map of Bahía Blanca. Four solutions are shown for the robust scenario, considering $\rho = 0.1$ and $\Gamma = 0.1$. Each solution displays the open collection points using different symbols for the two waste fractions: mixed and recyclable waste. The size of the symbol corresponds to the total number of installed bins (of any type), and the color of the symbol indicates the collection frequency selected for each collection point. The frequency value in the legend represents the number of days between consecutive visits.

In the solution with the minimum cost depicted in Figure 4(a), only a few collection points are open, with a low collection frequency. Moreover, all collection points, except two, are visited daily by the collection vehicle. Figure 4(b) and Figure 4(c) represent intermediate solutions that explore the trade-off between both goals. This trade-off is evident, on the one hand, in the larger number of open collection points and installed bins and, on the other hand, in the increment in the collection frequencies of open points, which have darker colors. Finally, the solution depicted in Figure 4(c) has the highest collection frequency and the highest cost for this scenario. This is evidenced by the large number of open collection points and installed bins, indicating an expensive solution, and the dark color with which most of the open points are represented, indicating a high collection frequency.

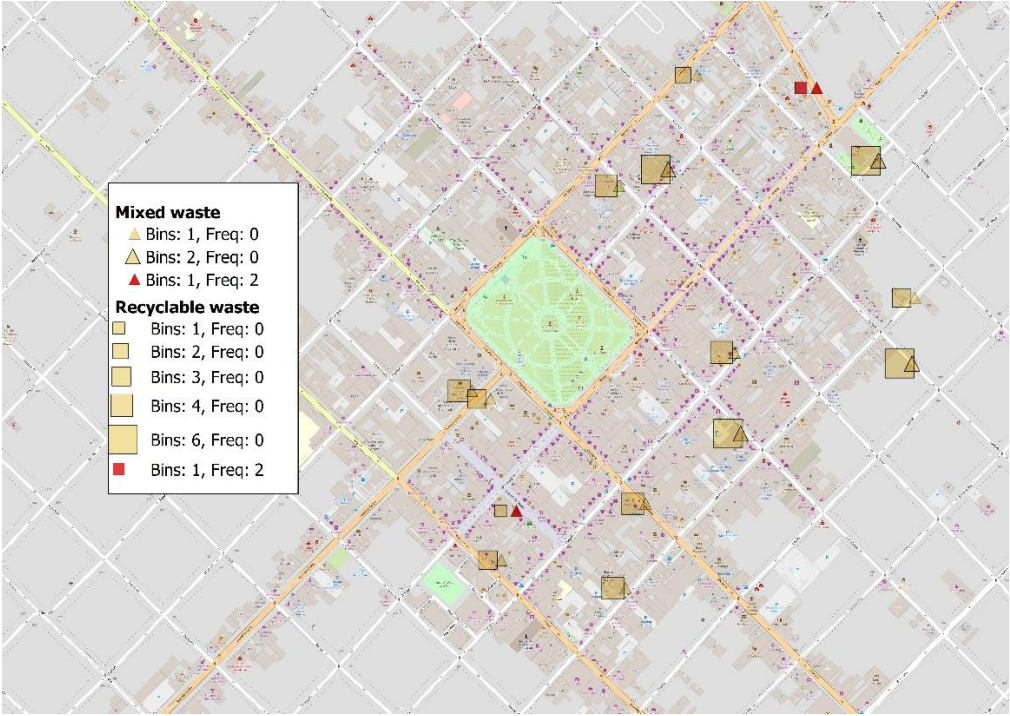

(a) Solution with minimum cost (with cost = 2342920 US$ and frequency = 27).

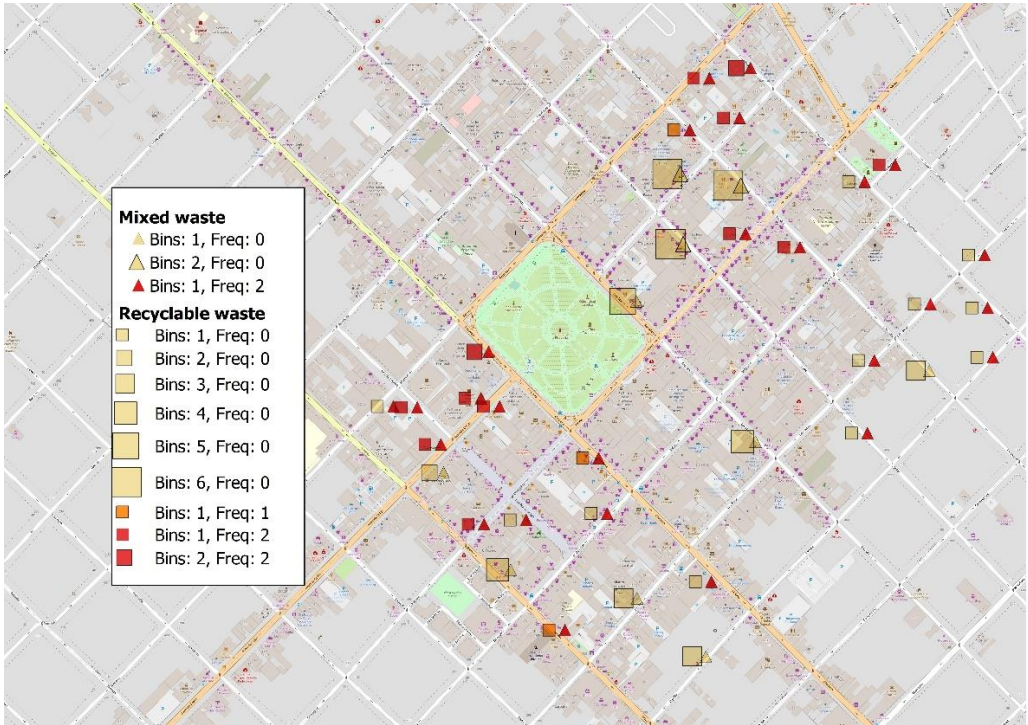

(b) Intermediate solution with cost = 3660000 US$ and frequency = 119.

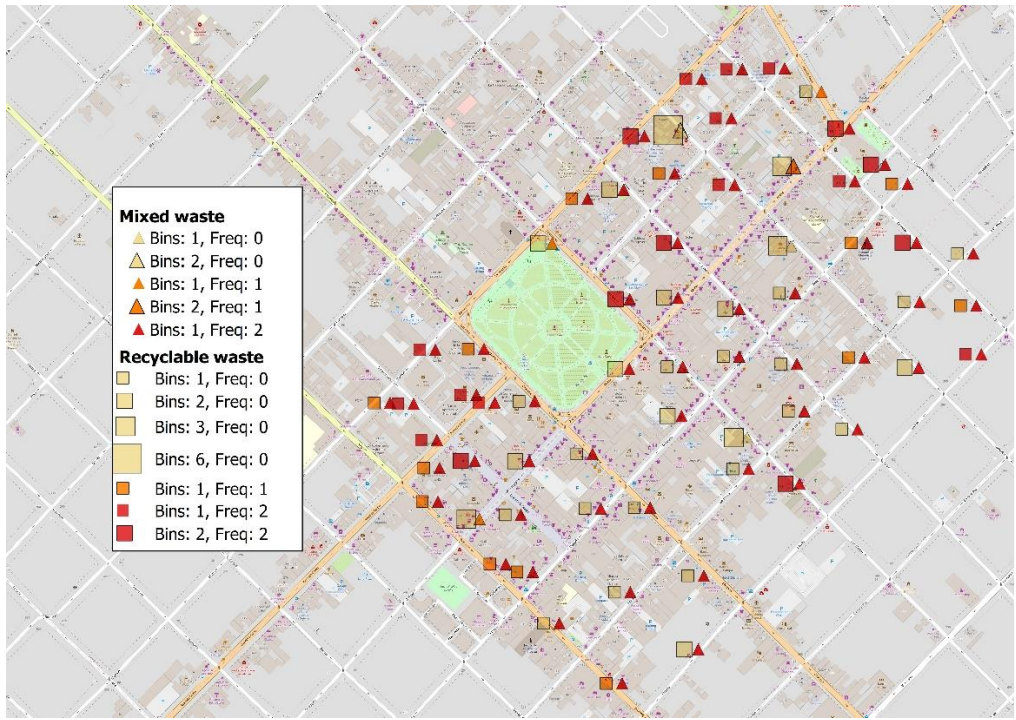

(c) Intermediate solution with cost = 5700000 US$ and frequency = 240.

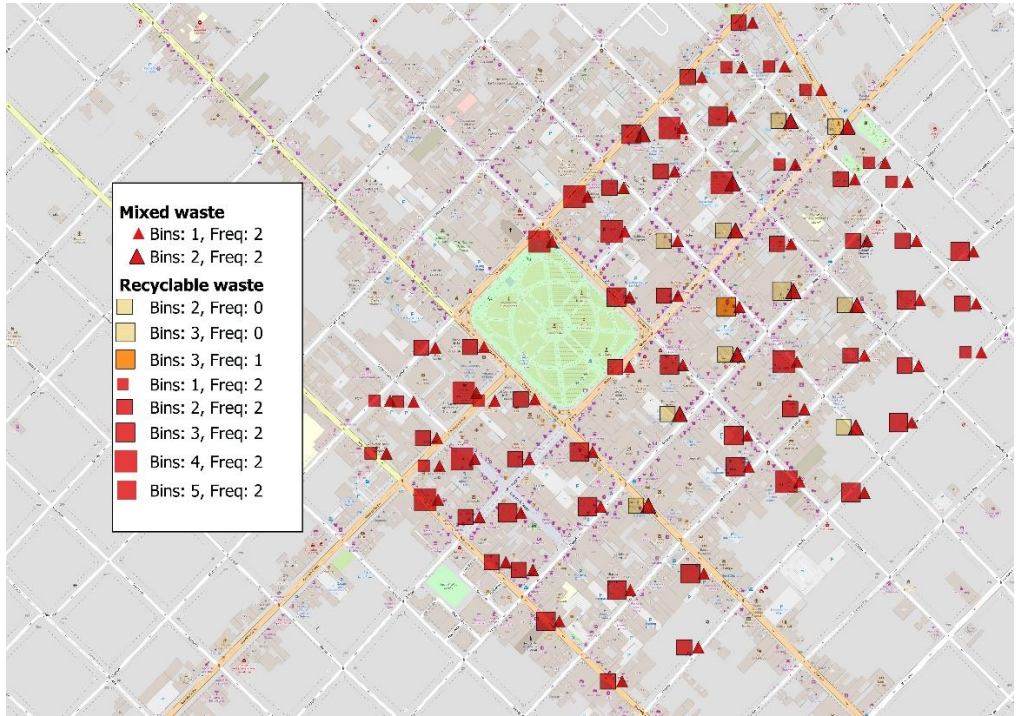

(d) Solution with maximum collection frequency (cost = 8380000 US$ and frequency = 331).

**Figure 4**: Geographic visualization of representative solutions of instance i73 considering $\rho = 0.1$ and $\Gamma = 0.1$. Source: background map OpenStreetMap[2].

## 4 Conclusions

Municipal solid waste management is a critical function of local governments who need to carefully manage the planning process to ensure a livable and sustainable urban environment. Effective planning and efficient operations, however, are not always straightforward due to a large amount of uncertainty surrounding key input parameters.

In an effort to address the uncertainty accompanying waste generation, this paper presents a robust mixed-integer mathematical programming model to help practitioners design the collection network – the first stage of an MSW system. In particular, a robust optimization model is proposed that aims to simultaneously minimize (i) the total cost of the network, including the establishment of collection points and the installment of waste bins, and (ii) the number of times that waste bins have to be emptied per week (or maximize the number of days that waste remains uncollected), which is a proxy for the routing costs. A threshold condition is imposed on the solution to ensure a minimum quality of service for users.

---

[2] https://www.openstreetmap.org/

Computational experiments are performed using a set of realistic instances based on a medium-size Argentinian city. The results obtained from the proposed model are compared with its deterministic counterpart for different scenarios of the robust model. The results show that, in general, the robust model is able to provide solutions that remain feasible with a high probability as the levels of waste generation vary throughout their practical range. Nevertheless, the robust solutions yield worse outcomes in terms of cost and collection frequency than those obtained from the deterministic model specially for scenarios with large waste generation variability. Moreover, the complexity of the computations, for what is already an NP-hard problem in its deterministic form, increases as the scenarios become more demanding, involving greater levels of uncertainty and variability.

Thus, this article contributes a robust model that allows practitioners to address the design of the collection network in an urban area while considering uncertainty in the amount of waste that is generated. Moreover, the robust model incorporates tuning parameters that allow practitioners to analyze the trade-off between the deterioration of optimization criteria compared to the deterministic model and the probability of violating the storage capacity in collection points due to uncertainty in waste generation rates. These parameters are based on the estimated variation in waste generation rates and the number of users considered to have varying waste generation rates. Overall, the main contributions of this work are twofold: it effectively addresses a real-world case study, providing insights that can be directly applied to practical waste management challenges, and it represents a pioneering effort in the related literature to address the computationally complex waste bin location problem considering an uncertain waste generation rate by means of robust optimization.

Future research might include approaching the waste management problem using a two- stage stochastic optimization model in which the distribution of waste generation is approximated by a given distribution. In that case, robust solutions obtained in this study could serve as an initial solution for the first stage of the stochastic model. Additionally, conducting simulation analyses to evaluate how well the proposed network of collection points in the robust solutions can accommodate different values of uncertain waste generation without requiring significant modifications would be valuable.

**Acknowledgements**

The work of the first author was supported by a Postdoctoral Scholarship co-funded by the Fulbright Foundation and the National Scientific and Technical Research Council of Argentina (CONICET).

**Appendix A. Robust transformation of the deterministic mixed-integer linear programming model.**

To introduce the concept of a robust solution in optimization problems, let's consider the following linear program:

$$\min_{x}\{c^T x | Ax \geq b\} \qquad (I)$$

In many optimization problems, the input data (A, b) is typically assumed to be known and certain. In this case, a feasible solution is defined as a solution $x'$ which satisfies the constraints $Ax' \geq b$.

However, in many real world applications input data is subjected to some degree of uncertainty[3]. In this regard, a robust feasible solution $x_r$ is a solution which remains feasible whatever realization of the input data within a reasonable prescribed "uncertainty set" $U$ in which the input data can vary (Ben-Tal and Nemirovski, 1999). Consequently, the robust counterpart to program (I) can be expressed as:

$$\min_x\{c^T x | Ax \geq b \; \forall (A,b) \in U\} \qquad (II)$$

Thus, $x_r$ is a feasible robust solution to program (II) if satisfies $Ax^r \geq b \; \forall \, (A,b) \in U$. Moreover, $x_r$ will be an optimal solution to the program (II) if $c^T x_r \leq c^T x'_r$ for any other robust feasible solution $x'_r$.

This definition of robust feasible solution can be extended to multi-objective linear optimization (Deb and Gupta, 2005; Kuroiwa and Lee, 2012) where instead of a scalar objective function we have a vector of objective functions $F(x)$. Let's consider the following multiobjective linear program:

$$\min_x\{F(x) | Ax \geq b \} \qquad (III)$$

where $F(x) = \{c_1^T x, c_2^T x, \dots, c_3^T x\}$. Similarly, to the single-objective program, the robust counterpart of program (III) can be expressed as:

$$\min_x\{F(x) | Ax \geq b \; \forall (A,b) \in U\} \qquad (IV)$$

Then, $x_r$ is a robust feasible solution to problem (IV) if it satisfies $Ax_r \geq b \; \forall (A,b) \in U$.

---

[3]If the vector $c$ is also uncertain, the equivalent formulation of (I): $\min_{x,t}\{y | c^T x \leq y, Ax \geq b\}$ can be considered. Without loss of generality, the uncertainty can be restricted to the constraints (Ben-Tal and Nemirovski, 1999).

Moreover, $x'_r$ is a robust efficient or non-dominated solution of problem (IV) if there is no other robust feasible solution $x_r$ of problem (IV) such that:

$$c^T x_r \leq c^T x'_r, i = 1, \ldots, n$$

$$\text{and } c^T x_r < c^T x'_r \text{ for some } 1 \leq j \leq n$$

In general, as in multiobjective deterministic problems (Deb and Deb, 2013), in multiobjective robust problems there is not a unique efficient solution (Deb and Gupta, 2005). The set of efficient solutions of a problem is called its "Pareto frontier".

Robust optimization techniques in order to find a robust feasible solution that remains feasible within the realizations of the input data in the uncertainty set impose a related cost on the optimal solution based on its degree of feasibility (Tirkolaee et al., 2020). Thus, despite being relatively more expensive, the derived robust solution is more reliable and will allow the decision-maker to have a practical solution for several realizations of the uncertain parameters which are difficult to predict. The first robust approach was proposed by Soyster (1973) using a linear optimization model that provides the best feasible solution for all possible realizations of random input data. This approach tends to find "over-conservative" solutions, which means that in order to ensure the robustness of the solutions in most realizations of the random input data the results are often suboptimal. To consider less conservative solutions, Ben-Tal and Nemirovski (1999, 2000) and El Ghaoui et al. (1998) proposed other approaches which involves solving the robust counterparts of the nominal problem using a quadratic objective function. Although these models can better approximate some types of uncertainties without defaulting to over-conservative solutions, they have the disadvantage of requiring the solution of a non-linear optimization problem, which tends to be more difficult to solve than the linear model of Soyster (1973). As an intermediate strategy to tackle both over-conservatism and high computational effort, Bertsimas and Sim (2004) introduced a methodology to control the conservatism level of the solutions with a linear formulation. We have adopted their methodology for our work.

As mentioned, one of the key parameters that can be affected by uncertainty is the waste generation rate (and composition of different waste fractions) (Singh, 2019). In the previous Model (1), the constraint that is affected if the waste generation rate of each fraction is random is (1e). In general, in studies that use robust optimization the uncertainty is controlled by two parameters: the uncertainty level and the conservatism level (Tirkolaee et al., 2020). In discussing the uncertainty level, initially defined by Ben-Tal and Nemirovski (2000), let's consider that the random waste generation rate $\tilde{b}_{im}$ varies symmetrically around its mean value $\bar{b}_{im}$ according to the following expression $\tilde{b}_{im} = (1 + \rho \xi_{im})\bar{b}_{im}$ where $\xi_{im}$ is a random variable distributed uniformly in the interval [-1,1] and $\rho$ is a positive constant. Then, $\rho$ is referred to as the uncertainty level as it quantifies the maximum variation of $\tilde{b}_{im}$. Consequently, $\tilde{b}_{im} \epsilon [\bar{b}_{im} - \rho \bar{b}_{im}, \bar{b}_{im} + \rho \bar{b}_{im}]$. For simplicity, we perform the following substitution $\hat{b}_{im} = \rho \bar{b}_{im}$ hereafter.

The other relevant parameter in the robust methodology proposed by Bertsimas and Sim (2004) is the conservatism level. The conservatism level is linked to the fact that it is unlikely that all the parameters that can potentially vary will actually do so simultaneously. Thus, the conservatism level fixes the number of parameters that are allowed to vary simultaneously. As aforementioned, in Model 1 the constraint that is affected by a random waste generation rate is (1e). For each Constraint (1e), let $J_{im}$ to be the set of waste generation rates that can vary (which will include the generator points that can be assigned to point $i$ due to the maximum distance constraint, i.e., $J_{im} = \{b_{jm}|j\epsilon I, di_{ji} \leq dis^{max}\}$). Then, for each Constraint (1e) we introduce a parameter $\Gamma_{im}$ that takes value within the interval $[0, |J_{im}|]$. This parameter $\Gamma_{im}$ is referred to as the conservatism level. Bertsimas and Sim (2004) proved that the robust solution obtained by applying their methodology will remain feasible as long as up to $\lfloor \Gamma_{im} \rfloor$ of the random coefficients change within their bounds, and up to one random coefficient changes up to $(\Gamma_{im} - \lfloor \Gamma_{im} \rfloor)\hat{b}_{jm}$. Moreover, Bertsimas and Sim (2004) demonstrated that, due to the symmetric distribution of variables, even if more than $\lfloor \Gamma_{im} \rfloor$ parameters vary, the robust solution will still be feasible with a very high probability. Furthermore, in their work Bertsimas and Sim showed that parameter $\Gamma_{im}$ controls the trade-off between the probability of violation of the constraints of the model due to uncertainty in the input parameters and the negative effect on the objective functions (compared to the value of the objective function in the deterministic problem). This is why parameter $\Gamma_{im}$ is also regarded as the "price of robustness". It represents the decay that must be accounted for in the objective function to attain a solution that remains feasible across multiple realizations of the uncertain parameters.

On this basis, to apply Bertsimas and Sim's robust methodology to Model (1), we replace deterministic Constraint (1e) with the following equations:

$$\sum_{i \in I} \bar{b}_{im} \left( \sum_{f \in F} Acc_f linV_{jimf} \right) \tag{A.1a}$$

$$+ \max_{\{S_{im} \cup t_{im} | S_{im} \subseteq J_{im}, |S_{im}| = \lfloor \Gamma_{im} \rfloor, t_{im} \in J_{im} \setminus S_{im}\}} \left\{ \sum_{j \in S_{im}} \hat{b}_{im} \left( \sum_{f \in F} Acc_f y_{jimf} \right) \right.$$

$$\left. + (\Gamma_{im} - \lfloor \Gamma_{im} \rfloor)\hat{b}_{t_{im}m} \left( \sum_{f \in F} Acc_f y_{t_{im}imf} \right) \right\} \leq \sum_{h \in H} cap_h v_{hmi}, \forall\, i \in I, m \in M$$

$$y_{jimf} \leq linV_{jimf} \leq y_{jimf} \tag{A.1b}$$

$$y \geq 0 \tag{A.1c}$$

where $y_{jimf}$ is defined as the variables that represents the absolute value of $linV_{jimf}$.

Two extreme cases can be used to better analyzed the proposed methodology. Note that when $\Gamma_{im} = 0, \forall\, i \in I, m \in M$ Constraint (A.1a) is reduced to Constraint (1e) of the deterministic Model

(1) since no parameter is allowed to vary. On the other hand, if $\Gamma_{im} = J_{im}, \forall\, i \in I, m \in M$ we are setting the number of parameters that can varies to the maximum, which results in the Soyster's method that generates "over-conservative" solutions as it is proved in Ben-Tal and Nemirovski (2000). Therefore, by varying $\Gamma_{im}$ in the interval $[0, |J_{im}|]$ we have the flexibility of adjusting the number of parameters that are allowed to vary.

According to Bertsimas and Sim, the maximization term of Constraint (A.1a) is called the "protection" term of the constraint since it maintains a gap between the storage capacity of the collection point (right hand side) and the deterministic accumulated waste term $\sum_{j \in I} \bar{b}_{jm}(\sum_{f \in F} Acc_f linV_{jimf})$. To linearize this protection term, Ben-Tal and Nemirovski proved that it can be replaced with the following linear program $\beta_{im}$ (A.2):

$$\beta_{im} = \max \sum_{j \in J_{im}} \left( \alpha_j \hat{b}_{jm} \left( \sum_{f \in F} Acc_f linV_{jimf} \right) \right) \qquad \text{(A.2a)}$$

Subject to

$$\sum_{j \in J_{im}} \alpha_j \leq \Gamma_{im} \qquad \text{(A.2b)}$$

$$\alpha_j \leq 1, \quad \forall j \in J_{im} \qquad \text{(A.2c)}$$

$$\alpha \geq 0 \qquad \text{(A.2d)}$$

To proof the equivalence of the protection term of Constraint (2a) and linear program $\beta_{im}$ (A.2), let's consider the optimal solution to $\beta_{im}$ (2). It can be deduced that the optimal solution is when $\lfloor \Gamma_{im} \rfloor$ of the $\alpha_j$ variables are equal to 1 and one $\alpha_j$ variable is equal to $\Gamma_{im} - \lfloor \Gamma_{im} \rfloor$. This is equivalent to selecting the subset $S_{im} \cup t_{im} | S_{im} \subseteq J_{im}, |S_{im}| = \lfloor \Gamma_{im} \rfloor, t_{im} \in J_{im} \setminus S_{im}$ which maximizes the previous protection term of Constraint (A.2a): $\sum_{j \epsilon S_{im}} \hat{b}_{im}(\sum_{f \epsilon F} Acc_f y_{jimf}) + (\Gamma_{im} - \lfloor \Gamma_{im} \rfloor)\hat{b}_{t_{im}m}(\sum_{f \epsilon F} Acc_f y_{t_{im}imf})$. Then, to continue with the linearization proposed by Ben-Tal and Nemirovski, we obtain the dual of program $\beta_{im}$ (A.2). Considering that $z$ and $p_i$ (with $i \epsilon J_{im}$) are the dual variables of Constraints (A.2b) and (A.2c), respectively, the dual model is defined as Model (A.3).

$$min \sum_{j \in J_{im}} p_{jm} + \Gamma z \qquad (A.3a)$$

Subject to

$$z + p_{jm} \geq \sum_{j \epsilon J_{im}} \left( \hat{b}_{jm} \left( \sum_{f \epsilon F} Acc_f y_{jimf} \right) \right), \quad \forall i \in J_{im} \qquad (A.3b)$$

$$z, p \geq 0 \qquad (A.3c)$$

By strong duality since $\beta_{im}$ (A.2) is feasible and bounded, then Model (A.3) is also feasible and bounded, and the optimal objective values of both models coincide. Thus, we can replace the protection term in Constraint (A.1a) with Model (A.3). Since Constraint (A.1a) is indexed in $i$ and $m$, variable $z$ and parameter $\Gamma$ are also indexed in these sets. On the other hand, variable $p$ is associated to each generation point that can deposit its waste in the specific potential collection point $i$. Finally, non-linear Constraint (A.1a) can be replaced with linear Constraints (A.4a)-(A.4c).

$$\sum_{i \epsilon I} \bar{b}_{im} \left( \sum_{f \epsilon F} Acc_f linV_{jimf} \right) + z_{im}\Gamma_{im} + \sum_{j \epsilon J_{im}} p_{jm} \leq \sum_{h \epsilon H} cap_h v_{hmi}, \quad \forall i \in I, m \in M \qquad (A.4a)$$

$$p_{jm} + z_{im} \geq \sum_{j \epsilon J_{im}} \left( \hat{b}_{jm} \left( \sum_{f \epsilon F} Acc_f y_{jimf} \right) \right), \quad \forall\, i,j \in I, m \in M \quad \text{(A.4b)}$$

$$z, p, y \geq 0 \quad \text{(A.4c)}$$

Thus, we consider the *Deterministic* Model as the model composed by Equations (1a)-(1l) and the *Robust* Model as the model composed by Equations (1a)-(1d), (1f)-(1l), (A.1b), and (A.4a)-(A.4c).

**Appendix B. Results of tests methods for obtaining approximation of the worst values of objectives within the Pareto frontier**

To apply the augmented ε-constraint method, it is necessary to determine the best and worst values of each objective function over the Pareto frontier (Mavrotas and Florios, 2013). However, obtaining the actual worst values of the objectives within the Pareto frontier can be challenging in integer programming multi-objective problems, and often approximations need to be used (Rossit et al., 2020). Different strategies were conducted to approximate these worst values. We employ lexicographic optimization as in Mavrotas and Florios (2013). In a bi-objective problem, the plain lexicographic approach involves optimizing the first criteria as a single-objective problem in the first stage. Then, in a second stage, the second criteria is optimized also with a single-objective optimization but including a constraint that prevents the first criteria from obtaining a worse value than the one achieved in the initial stage. In mixed-integer problems, the resolution of the second stage can be particularly challenging, as we are setting the first objective function to its optimal value (or near-optimal value if the initial stage has not converged) (Rossit et al., 2020). To enhance the lexicographic approach in the second stage, we used two different strategies: warm start and relaxation of the restriction on the first objective. The warm start strategy consists in the lexicographic approach using the solution of the first stage to initialize the optimization process of the second stage since the solution of the previous stage constitutes a feasible solution of the following stage. Providing the MIP solver with a feasible starting solution can greatly assist the solver by enabling efficient cuts in the branch and bound tree, effectively reducing the size of the problem to such an extent that further search in the branch and bound tree becomes possible (Pour et al., 2018). Then in the relaxation strategy, we relax the additional constraint that prevents the worsening of the first objective. This relaxation, as discussed in Rastegar and Khorram (2015), helps in obtaining feasible solutions by allowing more flexibility in the optimization process.

Results of the application of this enhancement strategies to the lexicographic approach are presented in Table A1 for the deterministic model and the robust model for a scenario with $\rho = 0.1$ and $\Gamma = 0.1$. Additionally, we studied the order of optimization of the objectives in the lexicographic approach. Thus, Table B.1 presents the results for each model considering both orders of optimization: cost in the first stage and frequency in the second, and vice versa. Starting from the left, the first column after the instance ID gives the cost result from the first run associated with optimizing the single objective function. The corresponding frequency value is given in the next column. The remaining columns report the results from the second run using four different strategies: straightforward lexicographic (L), lexicographic with relaxation (LR), lexicographic with relaxation and warm start (LRWS) and lexicographic with warm start (LWS). The results are expressed as a percent improvement over the results from the first run. A positive value indicates an improvement in the objective value compared to the first run, while a negative value implies a deterioration. Entries in bold denote the optimal solutions, i.e., when the solver returned a gap equal to zero.

**Table B.1.**: Results of the approximation of the worst values of objectives within the Pareto frontier.

| | | | | | | | | | | |
|---|---|---|---|---|---|---|---|---|---|---|
| | colspan="10" | **Deterministic model** |
| id | First run (opt. Cost) | | Second run (opt. Freq) | | | | | | | |
| | | | L | | LR | | LRWS | | LWS | |
| | Cost | Freq | Cost | Freq | Cost | Freq | Cost | Freq | Cost | Freq |
| i16 | **450000** | **6** | **0.00%** | **33.33%** | -0.65% | 33.33% | 0.00% | 33.33% | 0.00% | 33.33% |
| i58 | 2192920 | 20 | 0.00% | 10.00% | -4.30% | 50.00% | -4.24% | 50.00% | 0.00% | 10.00% |
| i73 | 2460000 | 22 | No solution | | -5.00% | 63.64% | -4.88% | 63.64% | 0.00% | 18.18% |
| i126 | 3665840 | 32 | No solution | | No solution | | -4.99% | 75.00% | 0.00% | 37.50% |
| i190 | 4412920 | 43 | No solution | | No solution | | -4.76% | 46.51% | 0.00% | 6.98% |
| id | First run (opt. Freq) | | Second run (opt. Cost) | | | | | | | |
| | | | L | | LR | | LRWS | | LWS | |
| | Cost | Freq | Cost | Freq | Cost | Freq | Cost | Freq | Cost | Freq |
| i16 | **4120600** | **96** | **54.13%** | **0.00%** | 57.05% | -4.17% | 57.05% | -4.17% | 54.13% | 0.00% |
| i58 | **12091140** | **348** | **26.81%** | **0.00%** | 33.51% | -4.89% | 33.51% | -4.89% | 26.81% | 0.00% |
| i73 | **14479980** | **438** | **27.69%** | **0.00%** | 34.12% | -4.79% | 34.12% | -4.79% | 27.69% | 0.00% |
| i126 | **19945400** | **756** | **17.12%** | **0.00%** | 24.19% | -4.89% | 24.19% | -4.89% | 17.12% | 0.00% |
| i190 | **44885920** | **1140** | **48.80%** | **0.00%** | 53.28% | -5.00% | 53.28% | -5.00% | 48.80% | 0.00% |
| | colspan="10" | **Robust model for $\rho = 0.1$ and $\Gamma = 0.1$** |
| id | First run (opt. Cost) | | Second run (opt. Freq) | | | | | | | |
| | | | L | | LR | | LRWS | | LWS | |
| | Cost | Freq | Cost | Freq | Cost | Freq | Cost | Freq | Cost | Freq |
| i16 | **450000** | **6** | **0.00%** | **33.33%** | -0.32% | 33.33% | 0.00% | 33.33% | 0.00% | 33.33% |
| i58 | **2231680** | **18** | No solution | | -4.85% | 50.00% | -4.98% | 61.11% | 0.00% | 11.11% |
| i73 | **2522920** | **22** | No solution | | -4.76% | 63.64% | -4.76% | 63.64% | 0.00% | 9.09% |
| i126 | **3815840** | **35** | No solution | | No solution | | -4.72% | 62.86% | 0.00% | 28.57% |
| i190 | **4488980** | **36** | No solution | | No solution | | -4.68% | 47.22% | 0.00% | 8.33% |
| id | First run (opt. Freq) | | Second run (opt. Cost) | | | | | | | |
| | | | L | | LR | | LRWS | | LWS | |
| | Cost | Freq | Cost | Freq | Cost | Freq | Cost | Freq | Cost | Freq |
| i16 | **450000** | **96** | 15.43% | 0.00% | 20.56% | -4.17% | 20.56% | -4.17% | **15.43%** | **0.00%** |
| i58 | **2231680** | **348** | 21.33% | 0.00% | 28.60% | -4.89% | 28.60% | -4.89% | **21.33%** | **0.00%** |
| i73 | **2522920** | **438** | 14.83% | 0.00% | 22.63% | -4.79% | 22.63% | -4.79% | **14.83%** | **0.00%** |
| i126 | **3815840** | **756** | 26.92% | 0.00% | 33.31% | -4.89% | 33.31% | -4.89% | **27.05%** | **0.00%** |
| i190 | **4488980** | **1140** | 27.16% | 0.00% | No solution | | 31.23% | -5.00% | **27.16%** | **0.00%** |

Results in Table B.1 show that, in general, the plain lexicographic (L) approach was able to obtain a solution within the threshold computing time when the objective of collection frequency is optimized first. However, when the cost is optimized in the first stage it only gets a solution in the smaller instances. The addition of relaxation (LR) exhibited a similar behaviour, performing better when the frequency is optimized in the first stage. On the other hand, the LRWS and LWS were able to obtain feasible solutions for every scenario, model and order of the objectives. In this regard, these two approaches demonstrated superior performance compared to the plain lexicographic approach and lexicographic approach with relaxation. Since LWS obtains feasible solutions without deteriorating the first objective, it is chosen as the preferred approach for the later results.